\documentclass[12pt]{amsart}
\usepackage{amssymb}
\usepackage{amsmath}

\theoremstyle{plain} 
\newtheorem{theorem}{Theorem}[section]
\newtheorem{corollary}[theorem]{Corollary}
\newtheorem{lemma}[theorem]{Lemma}
\newtheorem{proposition}[theorem]{Proposition}
\newtheorem{conjecture}[theorem]{Conjecture}

\theoremstyle{definition}

\newtheorem{example}[theorem]{Example}

\theoremstyle{remark}
\newtheorem{remark}[theorem]{Remark}

\numberwithin{equation}{section}

\def\Z{{\mathbb Z}}
\def\N{{\mathbb N}}
\def\R{{\mathbb R}}

\def\Q{{\mathbb Q}}
\def\I{{\mathbb I}}

\def\cB{{\mathcal B}}

\def\cD{{\mathcal D}}
\def\cF{{\mathcal F}}

\def\fkm{{\mathfrak m}}

\def\fkp{{\mathfrak p}}

\def\a{{\boldsymbol a}}
\def\b{{\boldsymbol b}}
\def\c{{\boldsymbol c}}
\def\d{{\boldsymbol d}}

\def\0{{\boldsymbol 0}}
\def\1{{\boldsymbol 1}}

\def\aalpha{{\boldsymbol{\alpha}}}
\def\bbeta{{\boldsymbol{\beta}}}

\def\llambda{{\boldsymbol{\lambda}}}

\def\codim{{\rm codim}}

\def\Gr{{\rm Gr}}

\def\Hom{{\rm Hom}}

\def\ZC{{\rm ZC}}

\title{Noetherian properties of rings of differential
operators of affine semigroup algebras
}

\author{Mutsumi Saito}
\address{Department of Mathematics,
Graduate School of Science,
Hokkaido University,
Sapporo, 060-0810, Japan
}
\email{saito@math.sci.hokudai.ac.jp}

\author{Ken Takahashi}
\address{Department of Mathematics,
Graduate School of Science,
Hokkaido University,
Sapporo, 060-0810, Japan}
\email{tomorrow.co.sky@hotmail.co.jp}

\date{January 19, 2007}

\begin{document}

\begin{abstract}
We consider the Noetherian properties of the ring of differential operators of an affine semigroup algebra.
First we show that it is always right Noetherian.
Next we give a condition, based on the data of the difference
between the semigroup and its scored closure, for 
the ring of differential operators being
anti-isomorphic to another ring of differential operators. 
Using this, we prove that the ring of differential operators is
left Noetherian if the condition is satisfied.
Moreover we give some other conditions 
for the ring of differential operators being 
left Noetherian.
Finally we conjecture necessary and sufficient conditions
for the ring of differential operators being 
left Noetherian.


\end{abstract}

\maketitle


\section{Introduction}
\label{INTRO}

Let $K$ be an algebraically closed field of characteristic zero.
Let $D(R)$ be the ring of differential operators of a finitely generated
commutative $K$-algebra $R$ as defined by Grothendieck \cite{ega4}.
We study the Noetherian properties of $D(R)$ when $R$ is an affine semigroup algebra.

It is well known that, if $R$ is a regular domain, then $D(R)$ is Noetherian,
 and the category of
left $D(R)$-modules and that of right $D(R)$-modules are
equivalent
(see for example \cite{Bjork}, \cite{McConnel-Robson}).
Bernstein-Gel'fand-Gel'fand \cite{BGG} showed that $D(R)$ is not
Noetherian in general if we do not assume the regularity.
However, $D(R)$ is known to be Noetherian 
for some families of interesting algebras;
Muhasky \cite{Muhasky} and Smith-Stafford \cite{Smith-Stafford}
independently proved that $D(R)$ is Noetherian if $R$ is an integral domain
of Krull dimension one. Tripp \cite{Tripp} proved that the ring $D(K[\Delta])$ of differential operators of the Stanley-Reisner ring $K[\Delta]$ is right Noetherian,
and gave a necessary and sufficient condition for $D(K[\Delta])$ to be
left Noetherian.

Let
$
A:=\{\, {\a}_1, {\a}_2,\ldots,
{\a}_n\,\}\subset \mathbb{Z}^d
$
be a finite subset.
We denote by $\mathbb{N} A$ the monoid generated by $A$, and
by $K[\mathbb{N} A]$ its semigroup algebra.
We consider
the ring $D(K[\N A])$ of differential operators of 
$K[\N A]$.
We saw in \cite{Prim,Saito-Traves} that the algebra
$D(K[\N A])$ is strongly related to
$A$-hypergeometric systems (also known as GKZ hypergeometric
systems), defined in \cite{KMM} and \cite{Hrabowski},
and systematically studied by Gel'fand and his collaborators
(e.g. \cite{Gelfand,GZK-def,GZK2,GKZ}).
The ring $D(K[\N A])$ was studied thoroughly when 
the affine semigroup algebra $K[\N A]$ is normal
(e.g. \cite{Jones} and \cite{Musson87,Musson94}).
In particular, Jones \cite{Jones} and Musson \cite{Musson94}
independently proved that  $D(K[\N A])$
is a Noetherian finitely generated $K$-algebra
when $K[\N A]$ is normal.
Traves and the first author \cite{Saito-Traves,Saito-Traves2} proved that
$D(K[\N A])$ is a finitely generated $K$-algebra
in general, and that
$D(K[\N A])$ is Noetherian if the semigroup $\N A$ is scored.
The scoredness means that $K[\N A]$ satisfies Serre's $(S_2)$ condition
and it is geometrically unibranched.
The question of Morita equivalence between $D(K[\N A])$ for a scored
$K[\N A]$ and that for its normalization
was studied in Smith-Stafford \cite{Smith-Stafford},
Chamarie-Stafford \cite{Chamarie-Stafford},
Hart-Smith \cite{Hart-Smith}, and
Ben Zvi-Nevins \cite{BenZvi-Nevins}.

Generally speaking, $D(R)$ is more apt to be right Noetherian
than to be left Noetherian, and a prime of height more than one is often
an obstacle for the left Noetherian property of $D(R)$
(see for example \cite{Cannings-Holland}, \cite{Musson-ZeroDivisors}, \cite{Smith-Stafford}, 
and \cite{Tripp}).
We observe this phenomenon as well.

As in \cite{Smith-Stafford}, by using Robson's lemma 
(Lemma \ref{Robson'sLemma}), we prove that
$D(K[\N A])$ is right Noetherian (Theorem \ref{theorem:RightNoetherian})
for any $A$.
To state conditions for the left Noetherain property, we need to introduce
the standard expression of a semigroup $\N A$;
let $S_c(\N A)$ be the {\it scored closure} of $\N A$, the smallest scored
semigroup containing $\N A$ (see \ref{def:S()}).
There exist $\b_1,\ldots,\b_m\in S_c(\N A)$ and faces $\tau_1,\ldots,\tau_m$ of $\R_{\geq 0}A$, the cone generated by $A$,
such that
\begin{equation}
\label{DefofStandardExpression}
\N A=S_c(\N A)\setminus \bigcup_{i=1}^m(\b_i+\Z(A\cap\tau_i)).
\end{equation}
Assuming that the expression \eqref{DefofStandardExpression} is irredundant, we see that 
the set $\{ \b_i+\Z(A\cap\tau_i)\,:\, i=1,\ldots, m\}$
is unique, and we call the expression \eqref{DefofStandardExpression}
the {\it standard expression} of $\N A$.

One way to prove the left Noetherian property is to show
the correspondence between left ideals and right ideals, and 
then to use the right Noetherian property.
To show the correspondence, we define a set $\cB$ 
based on the standard expression
(see \eqref{definition:cB} for the definition of $\cB$), 
and, when $\cB\neq\emptyset$, we consider a right $D(K[\N A])$-module
$K[\omega(\N A)]$, an analogue of the canonical module.
Then we see that the category of left $D(K[\N A])$-modules and that of
right $D(K[\omega(\N A)])$-modules are equivalent
(Theorem \ref{theorem:Left-Right}), and hence
we derive the left Noetherian property of $D(K[\N A])$ 
from the right Noetherian
property of $D(K[\omega(\N A)])$ (Theorem \ref{thm:LeftNoetherian}).
(For this reason and another technical reason, we prove the right Noetherian property not only of
$D(K[\N A])$ but of a little more general algebras.)
In this way, for example, we see that $D(K[\N A])$ is left Noetherian
if $K[\N A]$ satisfies Serre's $(S_2)$ condition.

Another way to prove the left Noetherian property is
the way similar to the one used for showing 
the right Noetherian property.
As in \cite[Proposition 7.3]{Smith-Stafford},
a sufficient condition for the left Noetherian property
is given in this way
(Theorem \ref{LeftNoetherianFacetAndOrigin}).

Finally we conjecture that $D(K[\N A])$ is left Noetherian
if and only if
for all $i$ with $\codim\, \tau_i>1$
$$
\left(
\bigcap_{\tau_j\succ\tau_i,\, \b_i-\b_j\in {K}\tau_j,\,
\codim\, \tau_j=1}\tau_j
\right)
=\tau_i.
$$
When this condition is not satisfied, we construct
a left ideal of $D(K[\N A])$ which is not finitely generated
(Theorem \ref{NecessaryCondition}).

This paper is organized as follows.
In Section \ref{Section:Definitions}, we recall some fundamental facts about
the rings of differential operators of semigroup algebras,
and fix some notation.
In Section \ref{section:FiniteGeneration},
we recall the results on the finite generation
in \cite{Saito-Traves,Saito-Traves2}, and generalize them
suitably for our proof of the right Noetherian property.
In Section \ref{section:PartialOrder}, we 
introduce preorders, which indicate
the $D$-module structures of $K[t_1^{\pm 1},
\ldots, t_d^{\pm 1}]$.
In Section \ref{section:RightNoether}, 
we prove the right Noetherian property.
In Section \ref{RightAndLeft},
we consider the correspondence between left $D$-modules and right $D$-modules.
In Section \ref{LeftNoether},
we consider the left Noetherian property;
we give a sufficient condition in Subsection \ref{subsection:SufficientCondition}
and a necessary condition
in Subsection \ref{subsection:NecessaryCondition}.

\section{Rings of Differential Operators}
\label{Section:Definitions}

In this section, 
we briefly recall some fundamental facts about
the rings of differential operators of semigroup algebras,
and fix some notation.

Let $K$ be an algebraically closed field of characteristic zero,
and $R$ a commutative $K$-algebra.
For $R$-modules $M$ and $N$, we inductively 
define the space of $K$-linear differential
operators from $M$ to $N$ of order at most $k$ by
$$
D^k(M, N):=\{ P\in \Hom_K(M,N)\,:\,
Pr-rP\in D^{k-1}(M,N)\,\,\text{for all $r\in R$}\}.
$$
Set $D(M,N):=\bigcup_{k=0}^\infty D^k(M, N)$,
and $D(M):=D(M,M)$.
Then, by the natural composition, $D(M)$ is a $K$-algebra, and
$D(M,N)$ is a $(D(N), D(M))$-bimodule.
We call $D(M)$ the {\it ring of differential operators}
of $M$.
For the generalities of the ring of differential operators,
see \cite{ega4}, \cite{McConnel-Robson}, \cite{Sw}, etc.

Let
\begin{equation}
A:=\{\, {\a}_1, {\a}_2,\ldots,
{\a}_n\,\}
\end{equation}
be a finite set of vectors in $\mathbb{Z}^d$.  Sometimes we identify
$A$ with the matrix of column vectors $(\a_1,\a_2,\ldots,\a_n)$.
Let $\N A$ and $\Z A$ denote the monoid and the group generated by
$A$, respectively.
{Throughout this paper, we assume that $\mathbb{Z} A=\mathbb{Z}^d$ for simplicity.}
{We also assume that the cone $\R_{\geq 0}A$ generated by $A$
is strongly convex.}

The (semi)group algebra of $\Z^d$ is the Laurent polynomial ring
$K[\mathbb{Z}^d]=K[t^{\pm 1}_1,\ldots,t^{\pm 1}_d]$.
Its ring of differential operators is
the ring of differential operators with Laurent polynomial coefficients
$$
D(K[\mathbb{Z}^d])=K[t^{\pm 1}_1,\ldots,t^{\pm 1}_d]\langle \partial_1,
\ldots,\partial_d\rangle,
$$
where $[\partial_i,t_j]=\delta_{ij}$,
$[\partial_i,t^{-1}_j]=-\delta_{ij}t^{-2}_j$, and the other
pairs of generators commute.
Here $[\,,\,]$ denotes the commutator, and $\delta_{ij}$ is 1 if $i=j$ and 0 otherwise.

The semigroup algebra ${K}[\mathbb{N} A]=\bigoplus_{\a\in \mathbb{N} A}{K} t^\a$ is
the ring of regular functions on the affine toric variety defined by
$A$, where $t^\a=t_1^{a_1}t_2^{a_2}\cdots t_d^{a_d}$ for
$\a={}^t(a_1,a_2,\ldots,a_d)$, the transpose of the row vector $(a_1,
\ldots, a_d)$.  
We say that $S\subseteq \Z^d$ is an {\it $\N A$-set}
if
${S}+\N A\subseteq {S}$.
Then ${K}[S]:=\bigoplus_{\a\in S}{K} t^\a$ is a $K[\N A]$-module.
Let $S, S'\subseteq \Z^d$ be $\N A$-sets.
Throughout this paper, we simply write $D(S,S')$
and $D(S)$ for $D(K[S], K[S'])$ and $D(K[S])$, respectively.
Then $D(S,S')$
can be
realized as a submodule 
of the ring $D(\mathbb{Z}^d)=D({K}[\mathbb{Z}^d])$ as follows:
$$
D(S, S')=\{
P\in D(\mathbb{Z}^d)\, :\, P({K}[S])\subset {K}[S']\}.
$$
(See \cite[p. 31]{Cannings-thesis}, 
\cite[Proposition 1.10]{Muhasky}, and
\cite[Lemma 2.7]{Smith-Stafford}.)

Put $s_j:=t_j\partial_j$ for $j=1,2,\ldots, d$.
We introduce a $\Z^d$-grading on the ring $D(\Z^d)$ as follows:
For $\a={}^t(a_1,a_2,\ldots,a_d)\in \mathbb{Z}^d$, set
$$
D(\Z^d)_\a:=
\{ P\in D(\Z^d)\, :\, [s_j, P]=a_jP\quad\mbox{for $j=1,2,\ldots,d$}\}.
$$
Then the algebra $D(\Z^d)$ is $\mathbb{Z}^d$-graded; $D(\Z^d)=\bigoplus_{\a\in \mathbb{Z}^d}D(\Z^d)_\a$.
Let $S, S'\subseteq \Z^d$ be $\N A$-sets.
For $\a\in \Z^d$, set $D(S, S')_\a:=D(S,S')\cap D(\Z^d)_\a$
and $D(S)_\a:=D(S)\cap D(\Z^d)_\a$.
Then $D(S)=\bigoplus_{\a\in \Z^d}D(S)_\a$ is a $\Z^d$-graded algebra,
and
$D(S, S')=\bigoplus_{\a\in \Z^d}D(S,S')_\a$ is a
$\Z^d$-graded $(D(S'), D(S))$-bimodule.
We can describe $D(S, S')_\a$ explicitly as in \cite[Theorem 2.3]{Musson87}.
For ${\a}\in \mathbb{Z}^d$, we define a subset
$\Omega_{S, S'}({\a})$ of $\Z^d$ by
$$
\Omega_{S,S'}({\a}) = \{\, {\d}
\in S: \; {\d} +\a \not\in S'\,\}
=
S\setminus (-{\a} +S').
$$
We simply write $\Omega_{S}({\a})$ for $\Omega_{S, S}({\a})$.
We regard the set $\Omega_{S, S'}({\a})$ as a set in $K^d$.

\begin{proposition}
\label{therem:Theorem 2.3Musson87}
$$D(S, S')_\a = 
t^{{\a}} 
{\mathbb{I}} (\Omega_{S,S'}({\a})),
$$
where
$$
\mathbb{I}(\Omega_{S,S'}({\a})):=\{
f(s)\in {K}[s]:={K}[s_1,\ldots,s_d]\,:\,
\mbox{$f$ vanishes on $\Omega_{S,S'}({\a})$}\}.
$$
\end{proposition}

\section{Finite Generation}
\label{section:FiniteGeneration}

In \cite{Saito-Traves} and \cite{Saito-Traves2},
we proved that $D(\N A)$ is finitely generated as a $K$-algebra,
and that $\Gr(D(\N A))$ is Noetherian if $\N A$
is scored, where $\Gr(D(\N A))$ is the graded ring
associated with the order filtration of $D(\N A)$.
In Section \ref{section:RightNoether}, we prove that
$D(\N A)$ is right Noetherian (Theorem \ref{theorem:RightNoetherian})
for any $A$,
by using Robson's lemma (Lemma \ref{Robson'sLemma}).
To this end, we need to generalize the results 
in
\cite{Saito-Traves} and \cite{Saito-Traves2}
in a wider situation.
This section is devoted to this purpose.

Let us recall the primitive integral support
function of a facet (face of codimension one) of the cone $\mathbb{R}_{\geq 0}A$.
We denote by $\mathcal{F}$ the set of facets of the cone
$\mathbb{R}_{\geq 0}A$.
Given $\sigma\in\mathcal{F}$,
we denote by $F_\sigma$ the {\it primitive integral support function\,}
of $\sigma$, i.e., $F_\sigma$ is a uniquely determined
linear form on $\mathbb{R}^d$ satisfying
\begin{enumerate}
\item
$F_\sigma(\mathbb{R}_{\geq 0}A)\geq 0$,
\item
$F_\sigma(\sigma)=0$,
\item
$F_\sigma(\mathbb{Z}^d)=\mathbb{Z}$.
\end{enumerate}

\begin{remark}
\label{remark:UpperBoundness}
Let $\sigma\in \cF$.
By the definition of $F_\sigma$, there exists $m\in \N$
such that $F_\sigma(\N A)\supseteq m+\N$.
Accordingly, for an $\N A$-set $S$, there exists $m\in \N$
such that $F_\sigma(S)\supseteq m+\N$.
\end{remark}

Let ${S_c}$ be a scored $\N A$-set.
Here we mean 
$$
{S_c}=\bigcap_{\sigma\in\cF}
\{ \a\in \Z^d\,:\, F_\sigma(\a)\in F_\sigma({S_c})\}.
$$
by `scored' .

\begin{lemma}
\label{lemma:Scored}
Let ${S_c}$ be a scored $\N A$-set.
Then
${S_c}$ is a finitely generated $\N A$-set if and only if
$F_\sigma({S_c})$ is a finitely generated $F_\sigma(\N A)$-set for each facet $\sigma$.
\end{lemma}

\begin{proof}
The only-if direction is obvious.

Suppose that $F_\sigma({S_c})$ is a finitely generated $F_\sigma(\N A)$-set for each facet $\sigma$.
Let
$M_\sigma:=\max (\N\setminus F_\sigma({S_c}))$
(cf. Remark \ref{remark:UpperBoundness}), and
$F_\sigma({S_c})_+:=\{ m\in F_\sigma({S_c})\,:\, m<M_\sigma\}\cup\{\infty\}$.
Then $F_\sigma({S_c})_+$ is a finite set.
For a map $\nu$ which assigns a facet $\sigma$ to an element of 
$F_\sigma({S_c})_+$, we define a subset ${S_c}(\nu)$ of ${S_c}$ by
$$
{S_c}(\nu)=\{ \a\in {S_c}\,:\, \text{$F_\sigma(\a)=\nu(\sigma)$
for all facets $\sigma$}\}.
$$
Here we agree that $F_\sigma(\a)=\infty$ means $F_\sigma(\a)\geq M_\sigma+1$.
Then
${S_c}=\bigcup_\nu {S_c}(\nu)$.
We also set
$$
{S_c}(\nu)_{\R}=\{ \a\in \R^d\,:\, \text{$F_\sigma(\a)=\nu(\sigma)$
for all facets $\sigma$}\}.
$$
Since ${S_c}$ is scored, ${S_c}(\nu)={S_c}(\nu)_{\R}\cap\Z^d$.

For each ray ($1$-dimensional face) $\rho$ of $\R_{\geq 0}A$,
fix $\d_\rho\in \N A\cap\rho$. 
Set
$$
F(\nu)=\{ \d_\rho\,:\, \text{$\rho\preceq \sigma$ for all facets $\sigma$ with
$\nu(\sigma)\neq\infty$}\}.
$$
Since any strongly convex cone is generated by its $1$-dimensional faces,
$$
\left\{ \a\in \R^d\,:\,
\begin{array}{ll}
F_\sigma(\a)=0 & \text{for all facets $\sigma$ with $\nu(\sigma)\neq\infty$}\\
F_\sigma(\a)\geq 0 & \text{for all facets $\sigma$ with $\nu(\sigma)=\infty$}
\end{array}
\right\}=
\R_{\geq 0}F(\nu).
$$
Hence $\R_{\geq 0}F(\nu)$ is the characteristic cone of ${S_c}(\nu)_{\R}$,
and there exists a polytope $P(\nu)$ such that
${S_c}(\nu)_{\R}=P(\nu)+\R_{\geq 0}F(\nu)$
(see \cite[\S 8.9]{Schrijver}).
Clearly there exists a finite set $G(\nu)$ such that
$\R_{\geq 0}F(\nu)\cap\Z^d=G(\nu)+\sum_{\d_\rho\in F(\nu)}\N \d_\rho$.
Thus ${S_c}(\nu)$ is generated by $(P(\nu)\cap\Z^d)\cup G(\nu)$ as an
$\N F(\nu)$-set.
Therefore ${S_c}$ is finitely generated as an $\N A$-set.
\end{proof}

Let ${S_c}$ be a scored finitely generated $\N A$-set.
Let $\b_1,\ldots,\b_m\in {S_c}$, and let $\tau_1,\ldots,\tau_m$ be 
faces of $\R_{\geq 0}A$.
Let
\begin{equation}
\label{DefOfS'}
{S}:={S_c}\setminus \bigcup_{i=1}^m(\b_i+\Z(A\cap\tau_i))
\end{equation}
satisfy $F_\sigma({S})=F_\sigma({S_c})$ for all facets $\sigma$.
We say that $S_c$ is the {\it scored closure\/} of $S$.
We assume that the expression \eqref{DefOfS'} is irredundant.
Then $\{ \b_i+\Z(A\cap\tau_i)\, :\, i=1,\ldots,m\}$ is unique,
and we call \eqref{DefOfS'} the {\it standard expression\/} of $S$.
We do not assume that ${S}$ is an $\N A$-set.
In the remainder of this section and the next section, 
$S_c$ and $S$ are fixed as above.

We define a ring of differential operators by
$$
D({S}):=\{
P\in D(\Z^d)\,:\,
P({K}[{S}])\subseteq {K}[{S}]\}.
$$

First we consider the $\Z^d$-graded structure of $D({S})$.
Put $D(S)_\a:=D(S)\cap D(\Z^d)_\a$.
Similarly to Proposition \ref{therem:Theorem 2.3Musson87}, we can write
$$
D({S})=\bigoplus_{\a\in\Z^d}D({S})_\a,\qquad
D({S})_\a=t^\a \I(\Omega_{{S}}(\a)),
$$
where
$\Omega_{{S}}(\a)={S}\setminus (-\a+{S})$.

\begin{proposition}
\label{prop:ZC1}
Let $\ZC$ stand for Zariski closure in ${K}^d$.
We regard $F_\sigma$ as a linear map from $K^d$ to $K$.
\begin{enumerate}
\item[{\rm (1)}]
$$
\ZC(\Omega_{S_c}(\d))=
\bigcup_{\sigma\in\cF}
\bigcup_{k\in F_\sigma({S_c})\setminus(-F_\sigma(\d)+F_\sigma({S_c}))}F_\sigma^{-1}(k).
$$
\item[{\rm (2)}]
$$
\ZC(\Omega_{{S}}(\d))=
\ZC(\Omega_{S_c}(\d))\cup
\bigcup_{\b_i-\d\in {S}+\Z(A\cap\tau_i)}
(\b_i-\d+{K}\tau_i).
$$
\end{enumerate}
In particular, $D({S})$ is a subalgebra of $D({S_c})$.
\end{proposition}

\begin{proof}
(1) is easy, and (2) follows 
from Lemma \ref{lem:3wayIntersection} below.
See also \cite[Proposition 5.1]{Saito-Traves2}.
\end{proof}

Fix $M\in \N$ so that
\begin{equation}
\label{Def-M}
M> \max F_\sigma({S_c})^c\cup\{ F_\sigma(\b_i)\, :\, i\}-
\min F_\sigma({S_c})\cup\{ F_\sigma(\b_i)\, :\, i\}
\quad\text{for all facets $\sigma$},
\end{equation}
where $F_\sigma({S_c})^c=\Z\setminus F_\sigma({S_c})$.
Then, for $N\geq M$,
$$
N+F_\sigma({S_c})\subseteq F_\sigma({S_c})
\quad\text{for all facets $\sigma$}.
$$

For a face $\tau$, put
$$
\overset{\circ\,\,\circ}{\N (A\cap\tau)}
:=
\{
\a\in \N (A\cap\tau)\, :\,
\text{$F_\sigma(\a)\geq M$ for all $\sigma\in\cF$
with $\sigma\not\succeq \tau$}
\}.
$$

\begin{lemma}
\label{lem:circ}
Let $\tau$ be a face of $\R_{\geq 0}A$. Then
$${S}+\overset{\circ\,\,\circ}{\N (A\cap\tau)}
\subseteq {S}.$$
\end{lemma}

\begin{proof}
Since ${S_c}$ is an $\N A$-set,
${S}+\N(A\cap\tau)\subseteq {S_c}$.
Let $\a\in \overset{\circ\,\,\circ}{\N (A\cap\tau)}$, and let $\b\in {S}$.
It remains to show that $\b+\a\notin \b_j+\Z(A\cap\tau_j)$ for any $j$.
If $\tau_j\not\succeq\tau$, then
there exists a facet $\sigma$ such that $\sigma\succeq \tau_j$ and
$\sigma\not\succeq \tau$.
For such $\sigma$,
$$
F_\sigma(\a+\b)
\geq
M+\min F_\sigma({S_c})> F_\sigma(\b_j).
$$
Hence $\a+\b\notin \b_j+\Z(A\cap\tau_j)$.

If $\tau_j\succeq\tau$, then
$\a+\b\notin \b_j+\Z(A\cap\tau_j)$ since
$\b\notin\b_j+\Z(A\cap\tau_j)$.
\end{proof}

\begin{lemma}
\label{lem:3wayIntersection}
$$
(\d+{S})\cap {S_c}\cap(\b_i+\Z(A\cap\tau_i))\neq\emptyset
\Leftrightarrow
\b_i-\d\in {S}+\Z(A\cap\tau_i).
$$
\end{lemma}

\begin{proof}
The implication $\Rightarrow$ is obvious.

For the implication $\Leftarrow$,
let $\a\in {S}$ and
$\d+\a\in\b_i+\Z(A\cap\tau_i)$.
For $\sigma\in\cF$
with $\sigma\succeq\tau_i$,
$$
F_\sigma(\d+\a)=F_\sigma(\b_i)\in F_\sigma({S_c}).
$$
If $\sigma\not\succeq\tau_i$, then
there exists $\a_\sigma\in \overset{\circ\,\,\circ}{\N (A\cap\tau_i)}$
such that
$F_\sigma(\d+\a+\a_\sigma)\in F_\sigma({S_c})$.
By Lemma \ref{lem:circ},
$\a+\sum_{\sigma\not\succeq\tau_i}\a_\sigma\in {S}$.
Hence
$\d+\a+\sum_{\sigma\not\succeq\tau_i}\a_\sigma$
belongs to the set in the left hand side.
\end{proof}

By Proposition \ref{prop:ZC1}, $\I(\Omega_{S_c}(\d))=\langle p_\d\rangle$,
where
\begin{equation}
\label{def:p_d}
p_\d(s)=
\prod_\sigma
\prod_{k\in F_\sigma({S_c})\setminus(-F_\sigma(\d)+F_\sigma({S_c}))}
(F_\sigma(s)-k).
\end{equation}

\begin{theorem}
\label{theorem:FiniteGeneration}
Let $\Gr(D({S_c}))$ denote the graded ring associated with the
order filtration of $D({S_c})$.
Then the following hold:
\begin{enumerate}
\item[{\rm (1)}]
$\Gr(D({S_c}))$ is finitely generated as a ${K}$-algebra.
$D({S_c})$ is left and right Noetherian.
\item[{\rm (2)}]
$D({S})$ is finitely generated as a ${K}$-algebra.
\item[{\rm (3)}]
$D({S_c})$ is finitely generated as a right $D({S})$-module.
\end{enumerate}
\end{theorem}

\begin{proof}
For (1) and (2),
the argument in \cite[Sections 5 and 6]{Saito-Traves2} works
with the new $M$ \eqref{Def-M}.
In \cite{Saito-Traves2}, we did the argument when ${S_c}=\R_{\geq 0}A\cap\Z^d$.

For (3),
we briefly recall some notation from \cite{Saito-Traves2}.
For a ray $\rho$ of the arrangement determined by $A$,
i.e. $\R\rho$ is the intersection of some hyperplanes $(F_\sigma=0)$,
take a nonzero vector ${\d}_\rho$
from $\Z^d\cap\rho$
satisfying the conditions:
\begin{equation}
\label{DefOfDRho-1}
\begin{array}{lll}
F_\sigma({\d}_\rho)\geq M\quad & \mbox{if} &
\quad F_\sigma({\d}_\rho)>0,\\
F_\sigma({\d}_\rho)\leq -M\quad &
\mbox{if} & \quad F_\sigma({\d}_\rho)<0,
\end{array}
\end{equation}
and
\begin{equation}
\label{DefOfDRho-2}
\d_\rho\in \Z (A\cap\tau)\cap\rho \quad
\mbox{for all faces $\tau$ of $\R_{\geq 0}A$ satisfying $\R\tau\supset\rho$.}
\end{equation}

Let $\mu$ be a map from $\mathcal{F}$ to
a set 
$$
\tilde{M}:=\{\, -\infty\,\}\cup\{\, +\infty\,\}
\cup \{\, m\in \Z\, :\, |m|< M\,\}.
$$
Define a subset $S_{c}(\mu)$ of $\Z^d$ by
\begin{equation}
S_{c}(\mu):=\{\,
{\d}\in \Z^d\, :\,
F_\sigma({\d})=\mu(\sigma)\quad
\mbox{for all $\sigma\in \mathcal{F}$}\,\},
\end{equation}
where we agree that $F_\sigma({\d})=
+\infty$ ($-\infty$, respectively) mean
$F_\sigma({\d})\geq M$ ($\leq -M$, respectively).
We also define
\begin{equation}
F(\mu)_{\R}:= \bigcap_{\sigma\in\cF} \left\{\,
{\d}\in \R^d\, :\,
\begin{array}{ll}
F_\sigma({\d})=0 &
\mbox{if $\mu(\sigma)\neq \pm \infty$}\\
F_\sigma({\d})\geq 0 &
\mbox{if $\mu(\sigma)= +\infty$}\\
F_\sigma({\d})\leq 0 &
\mbox{if $\mu(\sigma)= -\infty$}
\end{array}
\,\right\},
\end{equation}

In the case of (3), for $\d_1\in S_{c}(\mu)$ and $\rho\subset F(\mu)_\R$,
the deficiency ideal (cf. \cite[Definition 5.4]{Saito-Traves2})
is the same as in the case of (2),
i.e.,
$$
D({S_c})_{\d_1}\cdot D({S})_{\d_\rho}
=
D({S_c})_{\d_1+\d_\rho}\cdot
\I\left(\bigcup_{\b_i-\d_\rho\in {S}+\Z(A\cap\tau_i)}(\b_i-\d_\rho+{K}(A\cap\tau_i))\right).
$$
Hence the same argument in the proof of 
\cite[Theorem 5.14]{Saito-Traves2}
works as well.
\end{proof}

\section{Partial Preorders}
\label{section:PartialOrder}

In this section, we keep situation \eqref{DefOfS'},
and introduce preorders, which indicate
the $D(S)$-submodule structure of $K[\Z^d]=K[t_1^{\pm 1},
\ldots, t_d^{\pm 1}]$.

For an ideal $I$ of ${K}[s]$ and a vector $\c\in K^d$,
we define a new ideal $I+\c$ by
$$
I+\c:=
\{ f(s-\c)\, :\, f(s)\in I\}.
$$

The following lemma is immediate from the definition.

\begin{lemma}
\label{lemma4-1}
For a subset $V$ of $K^d$ and a point $\a\in K^d$, let $\I(V)$ denote the ideal of polynomials
vanishing on $V$, and $\fkm_\a:=\I(\{ \a\})$ the maximal ideal at $\a$.
Then the following hold:
\begin{enumerate}
\item[{\rm (1)}]
$\I(V)+\c=\I(V+\c)$.
\item[{\rm (2)}]
$\fkm_\a+(\b-\a)=\fkm_\b$.
\item[{\rm (3)}]
If $\fkp$ is prime, then so is $\fkp+\c$.
\end{enumerate}
\end{lemma}

Let $\fkp$ be a prime ideal of 
${K}[s]$.
In the set $\{ \fkp+\a\, :\,
\a\in \Z^d\}$, we define $\preceq_{{S}}$
by
\begin{equation}
\label{def:Preorder_S'}
\fkp\preceq_{{S}} \fkp+\a
\overset{\text{def.}}{\Leftrightarrow}
\I(\Omega_{{S}}(\a))\not\subseteq \fkp.
\end{equation}

\begin{remark}
When ${S}=\N A$, $\fkm_\a\preceq_{{S}} \fkm_\b$ if and only if
$\a\preceq\b$ in the sense of \cite{Saito-Traves}, which was also considered
in \cite{Musson-Van den Bergh}.
\end{remark}

\begin{lemma}
$\preceq_{{S}}$ is a partial preorder.
\end{lemma}

\begin{proof}
First, since $\I(\Omega_{{S}}(\0))=(1)\not\subseteq
\fkp$, we have $\fkp \preceq_{{S}} \fkp$.

Second, let $\fkp \preceq_{{S}} \fkp +\a$ and
$\fkp+\a \preceq_{{S}} \fkp +\a+\b$.
Then
we have $\I(\Omega_{{S}}(\a))\not\subseteq\fkp$
and $\I(\Omega_{{S}}(\b))\not\subseteq\fkp+\a$.
The latter is equivalent to
$\I(\Omega_{{S}}(\b)-\a)\not\subseteq\fkp$.
Since $\fkp$ is prime,
we have
$$
\I(\Omega(\a))\I(\Omega(\b)-\a)\not\subseteq\fkp.
$$

From the inclusion $D({S})_\b D({S})_\a\subseteq D({S})_{\a+\b}$,
we obtain
$$
\I(\Omega_{{S}}(\a))\I(\Omega_{{S}}(\b)-\a)\subseteq \I(\Omega_{{S}}(\a+\b)).
$$
Hence we have
$$
\I(\Omega_{{S}}(\a+\b))\not\subseteq\fkp,
$$
or equivalently
$$
\fkp\preceq_{{S}}\fkp +\a+\b.
$$
\end{proof}

\begin{lemma}
\label{lemma:poS'}
Let $\aalpha\in {K}^d$, and let $\a\in \Z^d$. Then
$\I(\tau)+\aalpha\preceq_{{S}} \I(\tau)+\aalpha+\a$
if and only if,
for all facets $\sigma\succeq\tau$,
$F_\sigma(\aalpha)\in F_\sigma({S_c})$
implies $F_\sigma(\aalpha+\a)\in 
F_\sigma({S_c})$, and,
for all faces $\tau_i\succeq\tau$,
$\aalpha+\a-\b_i\in {K}\tau_i$ implies
$\b_i-\a\notin {S}+\Z(A\cap\tau_i)$.
\end{lemma}

\begin{proof}
By definition,
$\I(\tau)+\aalpha\preceq_{{S}} \I(\tau)+\aalpha+\a$
means
$\ZC(\Omega_{{S}}(\a))\not\supseteq
\aalpha+{K}\tau$.
By Proposition \ref{prop:ZC1},
the latter condition means that,
for all facets $\sigma\succeq\tau$,
$F_\sigma(\aalpha)\in F_\sigma({S_c})$
implies $F_\sigma(\aalpha+\a)\in 
F_\sigma({S_c})$, and,
for all faces $\tau_i\succeq\tau$,
$\aalpha+\a-\b_i\in {K}\tau_i$ implies
$\b_i-\a\notin {S}+\Z(A\cap\tau_i)$.
\end{proof}

For $\aalpha\in {K}^d$ and a face $\tau$, set
\begin{equation}
\label{def:ES'tau}
E({S})_\tau(\aalpha):=
\{ \llambda\in {K}\tau/\Z(A\cap\tau)\,:\,
\aalpha-\llambda\in
{S}+\Z(A\cap\tau)\}.
\end{equation}
Define another partial preorder $\preceq_{{S},\tau}$ by
\begin{equation}
\label{def:PreorderS'tau}
\aalpha\preceq_{{S},\tau}
\bbeta\quad
\overset{\text{def.}}{\Leftrightarrow}\quad
E({S})_{\tau'}(\aalpha)\subseteq
E({S})_{\tau'}(\bbeta)
\quad \text{for all faces $\tau'$ with $\tau'\succeq\tau$}.
\end{equation}
We denote by $\aalpha\sim_{{S},\tau}\bbeta$
if $\aalpha\preceq_{{S},\tau}\bbeta$ and $\aalpha\succeq_{{S},\tau}\bbeta$,
or equivalently, if
$$
E({S})_{\tau'}(\aalpha)=
E({S})_{\tau'}(\bbeta)
\quad \text{for all faces $\tau'$ with $\tau'\succeq\tau$}.
$$

When ${S}=\N A$, the set $E({S})_\tau(\aalpha)$ was considered in 
\cite{IsoClass}.
As in the case when ${S}=\N A$, we have the following lemma.

\begin{lemma}
\label{lemma:FundamentalPropertyOfES'}
\begin{enumerate}
\item[{\rm (1)}]
$E({S})_\tau(\aalpha)$ is a finite set.
\item[{\rm (2)}]
$E({S})_{\R_{\geq 0}A}(\aalpha)=\{ \aalpha \mod \Z^d\}$.
\item[{\rm (3)}]
For a facet $\sigma\in \cF$,
$E({S})_\sigma(\aalpha)\neq\emptyset$
if and only if
$F_\sigma(\aalpha)\in F_\sigma({S})=F_\sigma({S_c})$.
\end{enumerate}
\end{lemma}

\begin{proof}
The proofs are the same as in the case when ${S}=\N A$.
See \cite[Propositions 2.2 and 2.3]{IsoClass}.
\end{proof}

By Lemma \ref{lemma:FundamentalPropertyOfES'} (2), $\aalpha+\Z^d=\bbeta+\Z^d$
if $\aalpha\preceq_{S,\tau}\bbeta$.

\begin{lemma}
\label{lemma:EqClassesAreFinite}
For any $\aalpha\in {K}^d$,
$\aalpha +\Z^d$ has only finitely many equivalence classes with respect to
$\sim_{{S},\tau}$.
\end{lemma}

\begin{proof}
Let $\bbeta\in \aalpha+\Z^d$.
If there exists $\llambda\in {K}\tau'$ such that
$\aalpha-\llambda\in \Z^d$, then
$\bbeta-\llambda\in \Z^d$, and
$E({S})_{\tau'}(\aalpha)$ and $E({S})_{\tau'}(\bbeta)$ are 
contained in the finite set $\llambda+\Q(\tau'\cap\Z^d)/\Z(A\cap\tau')$.
If there exists no such $\llambda\in {K}\tau'$,
then
$E({S})_{\tau'}(\aalpha)$ and $E({S})_{\tau'}(\bbeta)$ are empty.
Hence the number of equivalence classes is finite.
\end{proof}

Next we compare two preorders $\preceq_{{S}}$ and $\preceq_{{S},\tau}$.
By Lemmas \ref{lemma:poS'} and \ref{lemma:FundamentalPropertyOfES'},
\begin{eqnarray}
\label{conditionPreceqS1}
&& \I(\tau)+\aalpha\preceq_{{S}} \I(\tau)+\aalpha+\a\nonumber\\
&\Leftrightarrow&
\left\{
\begin{array}{l}
\text{for all facets $\sigma\succeq\tau$,
$E({S})_\sigma(\aalpha)\neq\emptyset
\Rightarrow E({S})_\sigma(\aalpha+\a)\neq\emptyset$}\\
\text{for all faces $\tau_i\succeq\tau$,
$\aalpha+\a-\b_i\notin E({S})_{\tau_i}(\aalpha)$.}
\end{array}
\right.
\end{eqnarray}
Note that
$\aalpha+\a-\b_i\notin E({S})_{\tau_i}(\aalpha+\a)$
is automatic.

Hence we have proved the following proposition.
\begin{proposition}
$$
\aalpha\preceq_{{S},\tau}\bbeta
\Rightarrow
\I(\tau)+\aalpha\preceq_{{S}}
\I(\tau)+\bbeta.
$$
\end{proposition}

We denote by $\I(\tau)+\aalpha\sim_{{S}}
\I(\tau)+\bbeta$
if
$\I(\tau)+\aalpha\preceq_{{S}}
\I(\tau)+\bbeta$ and $\I(\tau)+\aalpha\succeq_{{S}}
\I(\tau)+\bbeta$.

\begin{corollary}
\label{cor:equivalences}
$$
\aalpha\sim_{{S},\tau}\bbeta
\Rightarrow
\I(\tau)+\aalpha\sim_{{S}}
\I(\tau)+\bbeta.
$$
\end{corollary}

Similarly to \cite[Lemma 3.6]{Saito-Traves2},
the following holds.

\begin{lemma}
\label{lemma:S1Localization}
$$
{S}+\Z(A\cap\tau)
=[{S_c}+\Z(A\cap\tau)]
\setminus
\bigcup_{\tau_i\succeq \tau}
(\b_i+\Z(A\cap\tau_i)).
$$
\end{lemma}

\begin{lemma}
\label{lemma:S+Ztau}
\begin{eqnarray*}
{S_c}+\Z(A\cap\tau)
&=&
\{\a\in\Z^d \,:\,
F_\sigma(\a)\in F_\sigma({S})=F_\sigma({S_c})
\quad \text{for all facets $\sigma\succeq\tau$}
\}\\
&=&
\{\a\in\Z^d \,:\,
E({S})_\sigma(\a)\neq\emptyset\quad \text{for all facets $\sigma\succeq\tau$}
\}.
\end{eqnarray*}
\end{lemma}

\begin{proof}
This is immediate from the definitions and Lemma \ref{lemma:FundamentalPropertyOfES'}.
\end{proof}

\begin{theorem}
\label{theorem:Equivalence}
$$
\aalpha\sim_{{S},\tau}\bbeta
\Leftrightarrow
\I(\tau)+\aalpha\sim_{{S}}
\I(\tau)+\bbeta.
$$
\end{theorem}

\begin{proof}
It is left to prove the implication $\Leftarrow$.
We suppose that
$\I(\tau)+\aalpha\preceq_{{S}}
\I(\tau)+\aalpha+\a$.
Hence we suppose the two conditions
in \eqref{conditionPreceqS1}.
We show that $E({S})_{\tau'}(\aalpha)\subseteq E({S})_{\tau'}(\aalpha+\a)$
for all $\tau'\succeq\tau$.
We assume the contrary;
we suppose that
$\llambda\in E({S})_{\tau'}(\aalpha)\setminus E({S})_{\tau'}(\aalpha+\a)$.
Then we have
\begin{eqnarray}
&&\llambda\in{K}\tau'
\label{cond:a}\\
&&\aalpha-\llambda\in {S}+\Z(A\cap\tau')
\label{cond:b}\\
&&\aalpha+\a-\llambda\notin {S}+\Z(A\cap\tau').
\label{cond:c}
\end{eqnarray}
By \eqref{cond:a} and \eqref{cond:b},
$F_\sigma(\aalpha)\in F_\sigma({S_c})$ for all $\sigma\succeq\tau'$.
Then by \eqref{cond:a} and the first condition of \eqref{conditionPreceqS1}
$F_\sigma(\aalpha+\a-\llambda)=F_\sigma(\aalpha+\a)\in F_\sigma({S_c})$
for all $\sigma\succeq\tau'$.
Then by \eqref{cond:c} and Lemma \ref{lemma:S1Localization}
there exists $\tau_i\succeq\tau'$ such that
$\aalpha+\a-\llambda\in \b_i+\Z(A\cap\tau_i)$.
Hence
$\aalpha+\a-\b_i+\Z(A\cap\tau_i)=\llambda+\Z(A\cap\tau_i)
\in E({S})_{\tau_i}(\aalpha)$.
This contradicts the second condition of \eqref{conditionPreceqS1}.
\end{proof}

\section{Right Noetherian Property}
\label{section:RightNoether}

In this section, we assume that
$S_0$ is a scored finitely generated $\N A$-set, 
that the expression
$
S_m=S_0\setminus \bigcup_{i=1}^m(\b_i+\Z(A\cap\tau_i))
$
with all $\b_i\in S_0$ is irredundant, and that
$F_\sigma(S_m)=F_\sigma(S_0)$ for all facets $\sigma$.
In the notation in \eqref{DefOfS'}, $S_0=S_c$ and $S_m=S$.
We prove that $D(S_m)$ is right Noetherian by the induction on $m$.

Let
$$
D(S_0, S_m):=\{
P\in D(\Z^d)\,:\,
P({K}[S_0])\subseteq {K}[S_m]\},
$$
and for $\a\in \Z^d$ let
$$
D(S_0,S_m)_\a:=D(S_0,S_m)\cap D(\Z^d)_\a.
$$
Then
$D(S_0, S_m)$ is a right ideal of $D(S_0)$, and
a left ideal of $D(S_m)$.
We have
\begin{equation}
D(S_0, S_m)_\a=t^\a\I(\Omega_{S_0,S_m}(\a)),
\qquad
\Omega_{S_0,S_m}(\a)=S_0\setminus (-\a+S_m),
\end{equation}
and
\begin{equation}
\label{ZCfor0m}
\ZC(\Omega_{S_0, S_m}(\a))=
\ZC(\Omega_{S_0}(\a))\cup
\bigcup_{\b_i-\a\in S_0+\Z(A\cap\tau_i);\, 1\leq i\leq m}
(\b_i-\a+{K}\tau_i).
\end{equation}

We use the following Robson's lemma to prove the right Noetherian property
of $D(S_m)$.

\begin{lemma}[Proposition 2.3 in \cite{Robson}]
\label{Robson'sLemma}
Let $A$ be a right ideal of a right Noetherian ring $S$.
Let $R$ be a subring of $S$ containing $A$.
Suppose that $S$ is finitely generated as a right $R$-module,
and that
$S/A$ is a right Noetherian $R$-module.
Then
the ring $R$ is right Noetherian.
\end{lemma}

By Lemma \ref{Robson'sLemma} and Theorem \ref{theorem:FiniteGeneration},
we only need to show that
$D(S_0)/D(S_0, S_m)$ is a Noetherian right $D(S_m)$-module.

Let $k\leq m$, and let
$
S_k=S_0\setminus \bigcup_{i=1}^k(\b_i+\Z(A\cap\tau_i)).
$
Since we know, by Theorem \ref{theorem:FiniteGeneration},
 that $D(S_0)$ is right Noetherian,
and that
$D(S_0)$ is a finitely generated right $D(S_k)$-module,
the sequence of right ideals of $D(S_0)$
\begin{equation}
\label{TheSequence}
D(S_0, S_m)\subseteq D(S_0, S_{m-1})
\subseteq \cdots \subseteq
D(S_0, S_1)\subseteq D(S_0)
\end{equation}
is a sequence of finitely generated right $D(S_k)$-modules.

Set
\begin{equation}
\label{def:Mk}
M_k:=D(S_0, S_{k-1})/D(S_0, S_k).
\end{equation}
We want to show that each $M_k$ is a Noetherian
right $D(S_m)$-module.

\begin{lemma}
\label{lemma:Mk}
$$
\I(\Omega_{S_0, S_k}(\a))=
\langle p_\a\rangle \cdot
\bigcap_{\b_i-\a\in S_0+\Z(A\cap\tau_i);\, i\leq k}
\I(\b_i-\a+\tau_i).
$$
\end{lemma}

\begin{proof}
By \eqref{def:p_d} and \eqref{ZCfor0m},
$$
\I(\Omega_{S_0, S_k}(\a))=
\langle p_\a\rangle \cap
\bigcap_{\b_i-\a\in S_0+\Z(A\cap\tau_i);\, i\leq k}
\I(\b_i-\a+\tau_i).
$$

Suppose that $\b_i-\a\in S_0+\Z(A\cap\tau_i)$.
If $\b_i-\a+{K}\tau_i\subseteq
\ZC(\Omega_{S_0}(\a))$,
then there exists a facet $\sigma\succeq \tau_i$
such that $F_\sigma(\b_i-\a)\in F_\sigma(S_0)\setminus
(-F_\sigma(\a)+F_\sigma(S_0))$.
This contradicts the fact that $F_\sigma(\b_i)\in F_\sigma(S_0)$.
Hence $p_\a\notin \I(\b_i-\a+\tau_i)$.
If $fp_\a\in \cap_i \I(\b_i-\a+\tau_i)$,
then $f\in \cap_i \I(\b_i-\a+\tau_i)$ since $\I(\b_i-\a+\tau_i)$ is prime.
We have thus proved the assertion.
\end{proof}

\begin{corollary}
If $(M_k)_\a\neq 0$, then
$\b_k-\a\in S_0+\Z(A\cap\tau_k)$,
or equivalently,
if $(M_k)_{\b_k-\a}\neq 0$, then
$\a\in S_0+\Z(A\cap\tau_k)$.
\end{corollary}

\begin{proof}
This is immediate from Lemma \ref{lemma:Mk} and the definition of
$M_k$ \eqref{def:Mk}.
\end{proof}

\begin{lemma}
\label{lemma5.10}
If $(M_k)_{\b_k-\a}D(S_m)_\c\neq 0$,
then
$\a-\c\preceq_{S_m, \tau_k} \a$.
\end{lemma}

\begin{proof}
We have
\begin{eqnarray}
&&
\a-\c\not\preceq_{S_m,\tau_k} \a\nonumber\\
&\Leftrightarrow&
\I(\tau_k)+\a-\c\not\preceq_{S_m} \I(\tau_k)+\a\nonumber\\
&\Leftrightarrow&
\I(\Omega_{S_m}(\c))\subseteq
\I(\tau_k+\a-\c)\nonumber\\
&\Rightarrow&
(M_k)_{\b_k-\a}D(S_m)_\c=0.\label{LastImplication}
\end{eqnarray}
Here the first equivalence is by Theorem \ref{theorem:Equivalence};
the second is by the definition 
of $\preceq_{S_m}$ \eqref{def:Preorder_S'}.
For the implication \eqref{LastImplication},
let $X\in D(S_0,S_{k-1})_{\b_k-\a}$ and
$t^\c f(s)\in D(S_m)_\c$.
Since $M_k$ is a right $D(S_m)$-module,
$Xt^\c f(s)\in t^{\b_k-\a+\c}\I(\Omega_{S_0, S_{k-1}}(\b_k-\a+\c))$.
Since $f(s)\in \I(\Omega_{S_m}(\c))$, we have
$f(s)\in \I(\tau_k+\a-\c)$.
Hence by Lemma \ref{lemma:Mk}
$$
Xt^\c f(s)\in t^{\b_k-\a+\c}\I(\Omega_{S_0, S_{k-1}}(\b_k-\a+\c))
\cap\I(\tau_k+\a-\c)
\subseteq
t^{\b_k-\a+\c}\I(\Omega_{S_0, S_{k}}(\b_k-\a+\c)).
$$
Thus the implication \eqref{LastImplication} holds.
\end{proof}

The following proposition is immediate from Lemma \ref{lemma5.10}.

\begin{proposition}
\label{corollary5.11}
Let $C$ be a set of equivalence classes in $S_0+\Z(A\cap\tau_k)$
with respect to
$\sim_{S_m, \tau_k}$
such that $\d\preceq_{S_m, \tau_k}\c$ and $\c\in C$
imply $\d\in C$.

Then
$\bigoplus_{\c\in C}(M_k)_{\b_k-\c}$ is a right $D(S_m)$-submodule of
$M_k=\bigoplus_{\c\in S_0+\Z(A\cap\tau_k)}(M_k)_{\b_k-\c}$.
\end{proposition}

For $1\leq k\leq m$, set
\begin{equation}
\label{def:S-check}
\check{S_k}:=S_0\setminus \bigcup_{1\leq i\leq m;\, i\neq k}
(\b_i+\Z(A\cap\tau_i)).
\end{equation}
Then $D(\check{S_k})$ is a ${K}$-subalgebra of $D(S_0)$
(Prposition \ref{prop:ZC1}),
and we may assume that $D(\check{S_k})$ is right Noetherian by
the induction.
Since $M_k$ is a Noetherian right $D(S_0)$-module,
and $D(S_0)$ is finitely generated as a right $D(\check{S_k})$-module
by Theorem \ref{theorem:FiniteGeneration},
$M_k$ is a Noetherian right $D(\check{S_k})$-module.

The following lemma relates the Noetherian property as a right
$D(\check{S_k})$-module to that as a right $D(S_m)$-module.

\begin{lemma}
\label{Lemma5.13}
Let $C_1\subseteq C_2\subseteq \b_k-S_0+\Z(A\cap\tau_k)$,
and suppose that
$N_i:=\bigoplus_{\a\in C_i}(M_k)_\a$ $(i=1,2)$ are
right $D(S_m)$-submodules of $M_k$.
If, for any $\a, \a+\c \in C_2\setminus C_1$, $X\in (M_k)_\a$, and
$P\in D(\check{S_k})_\c$,
there exists $Q\in D(S_m)_\c$ such that $X.Q=X.P$,
then
$N_2/N_1$ is a Noetherian $\Z^d$-graded right $D(S_m)$-module.
\end{lemma}

\begin{proof}
Let $N$ be a $\Z^d$-graded right $D(S_m)$-submodule of $M_k$
with $N_1\subseteq N\subseteq N_2$.
Put
$$
\tilde{N}:=(\bigoplus_{\a\in C_2\setminus C_1}N_\a)D(\check{S_k}).
$$
Then $\tilde{N}$ is a right $D(\check{S_k})$-submodule of $M_k$.
By the assumption,
$\tilde{N}_\a=N_\a$ for all $\a\in C_2\setminus C_1$.
Hence $(\tilde{N}\cap N_2 + N_1)/N_1= N/N_1$.
Therefore
the Noetherian property of $M_k$ as a right $D(\check{S_k})$-module
implies that of $N_2/N_1$ as a $\Z^d$-graded right $D(S_m)$-module.
\end{proof}

Replacing $\a$ by $\b_k-\a$ in Lemma \ref{Lemma5.13},
we have the following corollary.

\begin{corollary}
\label{corollary5.14}
Let $C$ be an equivalence class in $S_0+\Z(A\cap\tau_k)$
with respect to $\sim_{S_m, \tau_k}$, and let $\d\in C$.
If, for any $\a, \a-\c \in C$, $X\in (M_k)_{\b_k-\a}$, and
$P\in D(\check{S_k})_\c$,
there exists $Q\in D(S_m)_\c$ such that $X.Q=X.P$,
then
$$
\bigoplus_{\a'\preceq_{S_m, \tau_k} \d}(M_k)_{\b_k-\a'}
/
\bigoplus_{\a'\prec_{S_m, \tau_k} \d}(M_k)_{\b_k-\a'}
$$
 is a Noetherian $\Z^d$-graded right $D(S_m)$-module,
where
$\prec_{S_m,\tau_k}$ means
$\preceq_{S_m,\tau_k}$
and
$\not\sim_{S_m,\tau_k}$.
\end{corollary}

\begin{proposition}
\label{proposition5.15}
For each equivalence class $C$ in $S_0+\Z(A\cap\tau_k)$
with respect to $\sim_{S_m, \tau_k}$,
the assumption in Corollary \ref{corollary5.14} is satisfied.
\end{proposition}

\begin{proof}
Similarly to Lemma \ref{lemma:Mk}, we see
\begin{eqnarray*}
\I(\Omega_{S_m}(\c))&=&
\langle p_\c\rangle \cdot
\I(\bigcup_{\b_i-\c\in S_m+\Z(A\cap\tau_i)}
(\b_i-\c+\tau_i))\\
\I(\Omega_{\check{S}_k}(\c))&=&
\langle p_\c\rangle \cdot
\I(\bigcup_{\b_i-\c\in \check{S}_k+\Z(A\cap\tau_i);\,
i\neq k}
(\b_i-\c+\tau_i)).
\end{eqnarray*}

Since $S_m\subseteq \check{S}_k$,
$\b_i-\c\in S_m+\Z(A\cap\tau_i)$ implies
$\b_i-\c\in \check{S}_k+\Z(A\cap\tau_i)$.
Hence, if $\b_k-\c\notin S_m+\Z(A\cap\tau_k)$,
then $\I(\Omega_{\check{S}_k}(\c))\subseteq 
\I(\Omega_{S_m}(\c))$, i.e.,
$D(\check{S_k})_\c\subseteq D(S_m)_\c$, and we have nothing to prove.

Suppose that $\b_k-\c\in S_m+\Z(A\cap\tau_k)$.
Let $f\in \I(\bigcup_{\b_i-\c\in \check{S}_k+\Z(A\cap\tau_i);\,
i\neq k}
(\b_i-\c+\tau_i))$.
If $f\in \I(\b_k-\c+\tau_k)$, then
$p_\c f\in \I(\Omega_{S_m}(\c))$, and again
we have nothing to prove.

Let $f\notin \I(\b_k-\c+\tau_k)$.
Let $X\in (M_k)_{\b_k-\a}$.
Suppose that $\b_k-\a\notin {K}\tau_k$.
Then there exists a facet $\sigma\succeq\tau_k$
such that $F_\sigma(\b_k-\a)\neq 0$.
Since $F_\sigma(s)-F_\sigma(\b_k-\c)\in
\I(\b_k-\c+\tau_k)$,
$p_\c f (F_\sigma(s)-F_\sigma(\b_k-\c))
\in \I(\Omega_{S_m}(\c))$.
We have
\begin{eqnarray*}
X.t^\c p_\c f (F_\sigma(s)-F_\sigma(\b_k-\c))
&=&
X.(F_\sigma(s-\c)-F_\sigma(\b_k-\c))t^\c p_\c f\\
&=&
X(F_\sigma(\a-\c)-F_\sigma(\b_k-\c))t^\c p_\c f\\
&=&
X F_\sigma(\b_k-\a)t^\c p_\c f
\end{eqnarray*}
Here the second equality above holds because
$$
X\in (M_k)_{\b_k-\a}
=t^{\b_k-\a} \I(\Omega_{S_0, S_{k-1}}(\b_k-\a))/t^{\b_k-\a}
\I(\Omega_{S_0, S_{k-1}}(\b_k-\a))
\cap\I(\a+\tau_k).
$$
Hence in this case
$$
X.t^\c p_\c f=
X.\frac{1}{F_\sigma(\b_k-\a)}
t^\c p_\c f (F_\sigma(s)-F_\sigma(\b_k-\c))
$$
as desired.

Finally suppose that $\b_k-\a\in {K}\tau_k$.
Since $\a\sim_{S_m, \tau_k}\a-\c$,
we have
$$
\a-(\a-\b_k)\in S_m+\Z(A\cap\tau_k)
\Leftrightarrow
\a-\c-(\a-\b_k)\in S_m+\Z(A\cap\tau_k),
$$
or equivalently,
$$
\b_k\in S_m+\Z(A\cap\tau_k)
\Leftrightarrow
\b_k-\c\in S_m+\Z(A\cap\tau_k).
$$
But the left hand side is false by the definition
of $(\b_k, \tau_k)$, and the right hand side
is true, which is one of our assumptions.
Hence the case when $\b_k-\a\in {K}\tau_k$
does not occur, and we have completed the proof of the
proposition.
\end{proof}

\begin{corollary}
\label{MkIsNoether}
$M_k$ is a Noetherian $\Z^d$-graded right $D(S_m)$-module.
\end{corollary}

\begin{proof}
Since $S_0+\Z(A\cap\tau_k)$ has only finitely many equivalence classes
with respect to $\sim_{S_m, \tau_k}$
by Lemmas \ref{lemma:EqClassesAreFinite} and \ref{lemma:S+Ztau}, $M_k$ is a
Noetherian $\Z^d$-graded right $D(S_m)$-module by Corollary \ref{corollary5.14}
and Proposition \ref{proposition5.15}.
\end{proof}

\begin{theorem}
\label{theorem:RightNoetherian}
$D(S_m)$ is right Noetherian.
\end{theorem}

\begin{proof}
By the sequence \eqref{TheSequence} and Corollary \ref{MkIsNoether},
$D(S_0)/D(S_0, S_m)$ is a Noetherian $\Z^d$-graded right $D(S_m)$-module, and hence $D(S_m)$ is $\Z^d$-graded right Noetherian by Robson's lemma.
Then, by the general theory of $\Z^d$-graded algebras (see \cite{Nastasescu-VanOystaeyen}),
$D(S_m)$ is right Noetherian.
\end{proof}
\section{Right Modules and Left Modules}
\label{RightAndLeft}

We retain the notation in Sections \ref{section:FiniteGeneration}
and \ref{section:PartialOrder}.
Thus ${S_c}$ is a finitely generated scored $\N A$-set,
we have an irredundant expression \eqref{DefOfS'}:
$${S}={S_c}\setminus \bigcup_{i=1}^m(\b_i+\Z(A\cap\tau_i))$$
with $\b_i\in {S_c}$,
and $F_\sigma({S_c})=F_\sigma({S})$ for all facets $\sigma$.
In this section, we assume that ${S}$ is an $\N A$-set.
When $S$ satisfies Serre's $(S_2)$ condition,
it is not difficult to see that $D(S)$ and $D(\omega(S))$ are
anti-isomorphic to each other, where 
$$
\omega(S) = -1\times(
\text{the weight set of the canonical module of $K[S]$}).
$$
Hence the left Noetherian property of $D(S)$ is derived from
the right Noetherian property of $D(\omega(S))$.
In this section, we give a sufficient condition for this argument
to stay valid.

For $P=\sum_\a t^\a f_\a(s)\in D(\Z^d)$,
the operator $P^\ast=\sum_\a f_\a(-s)t^\a$
is called the formal adjoint operator of $P$.
Then ${K}[\Z^d]={K}[t_1^{\pm 1},\ldots, t_d^{\pm d}]$
is a right $D({S})$-module by taking formal adjoint operators.

\begin{lemma}[cf. Proposition 4.1.5 in \cite{Saito-Traves}]
\label{Prop4.1.5}
Suppose that $\Lambda\subseteq \Z^d$ satisfies
that
$\a\in\Lambda$ and $\b\preceq_{{S},\, \{ \0\}}\a$
imply
$\b\in\Lambda$.
Then ${K}[ -\Lambda]$ is a right $D({S})$-submodule of
${K}[\Z^d]$.
\end{lemma}

\begin{proof}
Let $f_\a\in \I(\Omega_{{S}}(\a))$ and $\b\in \Lambda$.
Then
$(t^\a f_\a(s))^\ast .t^{-\b}
=f_{\a}(-s)t^\a.t^{-\b}
=f_\a(-s)t^{\a-\b}
=f_\a(\b-\a)t^{\a-\b}$.

The condition $\b-\a\not\preceq_{{S},\, \{ \0\}}\b$
 is equivalent to $\I(\Omega_{{S}}(\a))\subseteq
\fkm_{\b-\a}$ by Lemma \ref{lemma4-1} and Theorem \ref{theorem:Equivalence}.
Hence $(t^\a f_\a(s))^\ast .t^{-\b}=0$
if $\b-\a\not\preceq_{{S},\, \{ \0\}}\b$.
This proves the lemma.
\end{proof}

For the left Noetherian property, we construct an $\N A$-set
$\omega({S})$, and show a duality between $D({S})$ and $D(\omega({S}))$.
To construct $\omega({S})$, we prepare some notation.
Let $\tilde{\cF}$ denote the union
$\cF\cup\{\tau_i\,:\, i=1,\ldots,m\}$.
Set
\begin{equation}
\label{definition:cB}
\cB:=
\left\{
(\b_\tau)_{\tau\in\tilde{\cF}}\,:\,
\begin{array}{l}
\bullet\quad \text{$\b_\tau\in {K}\tau\cap\Z^d/\Z(A\cap\tau)$
for all $\tau\in\tilde{\cF}$.}\\
\bullet\quad
\text{For all $i$ and all $\tau\in\tilde{\cF}$ with $\tau\succeq\tau_i$,
there exists $j$ with}\\
\quad\text{$\,\,\,\tau_j=\tau$ such that
$\b_i+\b_{\tau_i}=\b_j+\b_{\tau_j}$ mod $\Z(A\cap\tau)$.}
\end{array}
\right\}.
\end{equation}

Throughout this section, we assume 
\begin{equation}
\label{condition:LeftNoether}
\cB\neq\emptyset.
\end{equation}

Fix an element $(\b_\tau)\in\cB$ once for all.
We define a subset $\omega({S})$ of $\Z^d$ by
\begin{eqnarray*}
\omega({S})&:=&
\{ \a\in \Z^d\,:\,
\text{$\b_\tau\notin E({S})_\tau(-\a)$ for each $\tau\in\tilde{\cF}$}
\}\\
&=&
\{ \a\in \Z^d\,:\,
\text{$-\a-\b_\tau\notin {S}+\Z(A\cap\tau)$ for each $\tau\in\tilde{\cF}$}
\}.
\end{eqnarray*}

By Lemma \ref{Prop4.1.5},
${K}[\omega({S})]$ is a right $D({S})$-module.

\begin{remark}
If ${S}$ satisfies Serre's $(S_2)$ condition,
\begin{equation}
\label{S_2}
{S}=\bigcap_{\sigma\in\cF}({S}+\Z(A\cap\sigma)),
\end{equation}
then $\tilde{\cF}=\cF$, and $(\0)\in\cB$.
Hence condition \eqref{condition:LeftNoether} is satisfied.

When $\N A$ satisfies Serre's $(S_2)$ condition,
$-\omega(\N A)$ for $(\0)\in\cB$ is the weight set
of a right $D(\N A)$-module $H_\fkm^d({K}[\N A])^\ast$,
the Matlis dual of the local cohomology module $H_\fkm^d({K}[\N A])$.
\end{remark}

\begin{lemma}
\label{omegaS'}
$\omega({S})$ is an $\N A$-set.
\end{lemma}

\begin{proof}
Let $\a\in\omega({S})$, and $\b\in \N A$.
Suppose that $\a+\b\notin \omega({S})$.
Then there exists a face $\tau\in \tilde{\cF}$
such that $-\a-\b-\b_\tau\in {S}+\Z(A\cap\tau)$.
Then from the $\N A$-stability of ${S}$ 
we obtain $-\a-\b_\tau\in {S}+\Z(A\cap\tau)$,
which contradicts the assumption $\a\in\omega({S})$.
\end{proof}

\begin{lemma}
\label{lemma:AdjointOperators}
$D({S})^\ast\subseteq D(\omega({S}))$,
where
$$
D({S})^\ast=\{
P^\ast\,:\,
P\in D({S})\}.
$$
\end{lemma}

\begin{proof}
Since the $\N A$-set ${S_c}$ is finitely generated,
$\omega({S})$ is not empty.
We know that $\omega({S})$ is an $\N A$-set by Lemma \ref{omegaS'}.
Hence, if $P^\ast\in D({S})^\ast$ satisfies
$P^\ast(t^\a)=0$ for all $\a\in\omega({S})$,
then $P^\ast=0$.
\end{proof}

Next we show that $\omega({S})$
is of the form considered in Sections \ref{section:FiniteGeneration}
and \ref{section:PartialOrder}.
Then we show that $D({S})^\ast= D(\omega({S}))$
under condition \eqref{condition:LeftNoether}.
Thus we deduce the left Noetherian property
of $D({S})$ from the right Noetherian property of
$D(\omega({S}))$ if condition \eqref{condition:LeftNoether} is satisfied.

We define the scored closure $S_c(\omega({S}))$ of $\omega({S})$ by
\begin{equation}
\label{def:S()}
S_c(\omega({S})):=
\bigcap_{\sigma\in\cF} \{ \a\in \Z^d\,:\,
F_\sigma(\a)\in F_\sigma(\omega({S}))\}.
\end{equation}

\begin{lemma}
\label{lemma:FiniteGenerationOfS()}
$S_c(\omega({S}))$ is a finitely generated scored $\N A$-set.
\end{lemma}

\begin{proof}
Clearly $S_c(\omega({S}))$ is scored.
Since $\omega({S})$ is an $\N A$-set by Lemma \ref{omegaS'},
$S_c(\omega({S}))$ is also an $\N A$-set.

For the finite generation, 
by Lemma \ref{lemma:Scored}, it is enough to prove that
each $F_\sigma(S_c(\omega({S})))=F_\sigma(\omega({S}))$ is finitely generated as an $F_\sigma(\N A)$-set.
For this, it suffices to show that $F_\sigma(\omega({S}))$ is bounded below.
For each facet $\sigma$, we know by Lemma \ref{lemma:S1Localization}
$$
{S}+\Z(A\cap\sigma)=
[{S_c}+\Z(A\cap\sigma)]\setminus
\bigcup_{\tau_i=\sigma}(\b_i+\Z(A\cap\sigma)).
$$
Clearly
$$
F_\sigma(\omega({S}))\subseteq
-(F_\sigma({S})^c\cup \{ F_\sigma(\b_i)\,:\, \tau_i=\sigma\}).
$$
This proves the finite generation, since the right hand side is bounded below.
\end{proof}

\begin{corollary}
\label{cor:omegaS'isRightNoetherian}
$D(\omega({S}))$ is right Noetherian.
\end{corollary}

\begin{proof}
By the Noetherian property of ${K}[\N A]$,
Lemma \ref{lemma:FiniteGenerationOfS()},
and the similar argument to that in 
\cite[Proposition 3.4]{Saito-Traves2},
$\omega({S})$
can be written of the form considered 
in Sections \ref{section:FiniteGeneration}
and \ref{section:PartialOrder}.
Hence we can apply
Theorem \ref{theorem:RightNoetherian} to $\omega({S})$.
\end{proof}

\begin{lemma}
\label{lemma:LocalizationAlongFacet}
For $\tau\in\tilde{\cF}$,
$$
\omega({S})+\Z(A\cap\tau)
=
\{
\a\in\Z^d\,:\,
\text{
$E({S})_{\tau'}(-\a)\not\ni \b_{\tau'}$
for any $\tau'\in\tilde{\cF}$ with $\tau'\succeq\tau$}
\}.
$$
\end{lemma}

\begin{proof}
The inclusion `$\subseteq$' is clear by definition.

For the inclusion `$\supseteq$',
let $\a\in\Z^d$ satisfy $E({S})_{\tau'}(-\a)\not\ni \b_{\tau'}$
for any $\tau'\in\tilde{\cF}$ with $\tau'\succeq\tau$.
Then we can take $\b\in \N(A\cap\tau)$ such that
$F_{\sigma}(-\a-\b)\notin F_{\sigma}({S})$ for any facet 
$\sigma\not\succeq\tau$.
Then $E({S})_{\tau'}(-\a-\b)=\emptyset$ for any face $\tau'\not\succeq\tau$.
In particular, $\a+\b\in \omega({S})$.
\end{proof}

\begin{lemma}
\label{lemma:FacetE-sets}
Let $\sigma$ be a facet.
For $\llambda\in\Z^d\cap{K}\sigma$,
$$
\llambda\in E(\omega({S}))_\sigma(-\a)
\Leftrightarrow
\b_\sigma-\llambda
\notin
E({S})_\sigma(\a).
$$
\end{lemma}

\begin{proof}
By Lemma \ref{lemma:LocalizationAlongFacet},
$$
\omega({S})+\Z(A\cap\sigma)=
\{ \a\in \Z^d\,:\,
-\a\notin \b_\sigma+{S}+\Z(A\cap\sigma)\}.
$$
Hence
$$
\Z^d
=
[\omega({S})+\Z(A\cap\sigma)]
\coprod -[\b_\sigma+{S}+\Z(A\cap\sigma)].
$$
Hence the assertion follows.
\end{proof}

The following proposition may be considered as the duality between
${S}$ and $\omega({S})$.

\begin{proposition}
\label{ReflectionOfS}
\begin{eqnarray*}
{S}
&=&
\{
\a\in\Z^d\,:\,
\text{$E({S})_{\tau}(\a)\ni\0$ for all $\tau\in\tilde{\cF}$}
\}\\
&=&
\{
\a\in\Z^d\,:\,
\text{$E(\omega({S}))_{\tau}(-\a)\not\ni\b_\tau$ for any $\tau\in\tilde{\cF}$}
\}\\
&=& \omega(\omega({S})).
\end{eqnarray*}
\end{proposition}

\begin{proof}
For any face $\tau$,
we have
$$
{S}+\Z(A\cap\tau)=
{S_c}+\Z(A\cap\tau)\setminus
\bigcup_{\tau_i\succeq\tau}(\b_i+\Z(A\cap\tau_i)).
$$
Hence
$$
{S}=
S_2({S})\cap
\bigcap_{\codim\tau_i>1}({S}+\Z(A\cap\tau_i)),
$$
where
$S_2({S})=\bigcap_{\sigma\in\cF}({S}+\Z(A\cap\sigma))$
is the $S_2$-closure of ${S}$.
This means the first equality of the proposition.

For the second equality, first note that
by Lemma \ref{lemma:FacetE-sets}
\begin{equation}
E({S})_\sigma(\a)\ni\0
\Leftrightarrow
E(\omega({S}))_\sigma(-\a)\not\ni\b_\sigma
\end{equation}
for any facet $\sigma$.
We have
\begin{eqnarray*}
S_2({S})
&=&
\{\a\in\Z^d\,:\,
E({S})_\sigma(\a)\ni\0\text{\quad for all facets $\sigma$}\}\\
&=&
\{\a\in\Z^d\,:\,
E(\omega({S}))_\sigma(-\a)\not\ni\b_\sigma\text{\quad for any facet $\sigma$}\}.
\end{eqnarray*}

Suppose that $\a\in S_2({S})$ and
$\b_{\tau_i}\in E(\omega({S}))_{\tau_i}(-\a)$ for some $\tau_i$
with $\codim\tau_i>1$.
Then $-\a-\b_{\tau_i}\in \omega({S})+\Z(A\cap\tau_i)$.
By Lemma \ref{lemma:LocalizationAlongFacet},
$E({S})_{\tau_i}(\a+\b_{\tau_i})\not\ni\b_{\tau_i}$,
or equivalently,
$E({S})_{\tau_i}(\a)\not\ni\0$.
We have thus proved the inclusion `$\subseteq$' of the second equation.

Since ${S}=S_2({S})\setminus\bigcup_{\codim\tau_i>1}(\b_i+\Z(A\cap\tau_i))$,
and since the right hand side of the second equality is included in
$S_2({S})$, to prove the inclusion `$\supseteq$',
it suffices to show that
$\b\in\b_i+\Z(A\cap\tau_i)$ with $\codim\tau_i>1$ does not belong to
the right hand side.
Since $(\b_\tau)$ belongs to $\cB$,
for any $\tau\in\tilde{\cF}$ with
$\tau\succeq\tau_i$ there exists $j$ with $\tau_j=\tau$
such that
$$
\b+\b_{\tau_i}-\b_\tau\in \b_j+\Z(A\cap\tau).
$$
In particular,
$$
\b+\b_{\tau_i}-\b_\tau\not\in {S}+\Z(A\cap\tau).
$$
Hence
$$
\b_\tau\notin E({S})_\tau(\b+\b_{\tau_i})
\qquad
\text{for any $\tau\in\tilde{\cF}$ with $\tau\succeq\tau_i$.}
$$
This means by Lemma \ref{lemma:LocalizationAlongFacet}
$$
-\b-\b_{\tau_i}\in\omega({S})+\Z(A\cap\tau_i).
$$
This is equivalent to
$$
\b_{\tau_i}\in E(\omega({S}))_{\tau_i}(-\b).
$$
Hence $\b$ does not belong to the right hand side of the second equality
of the proposition.
\end{proof}

\begin{theorem}
\label{theorem:SurjectivityByS_2}
Under condition \eqref{condition:LeftNoether}, 
$$
D(\omega({S}))=D({S})^\ast.
$$
\end{theorem}

\begin{proof}
By Lemma \ref{lemma:AdjointOperators} and Proposition \ref{ReflectionOfS},
$D(\omega({S}))^\ast\subseteq D({S})$.
Hence
$$
D(\omega({S}))=
D(\omega({S}))^{\ast\ast}
\subseteq
D({S})^\ast
\subseteq
D(\omega({S})).
$$
Hence $D(\omega({S}))=D({S})^\ast$
and
$D({S})=D(\omega({S}))^\ast$.
\end{proof}

\begin{theorem}
\label{theorem:Left-Right}
Assume that ${S}$ satisfies condition \eqref{condition:LeftNoether}.
Then there exist one-to-one correspondences between
left modules, left ideals, right modules, right ideals of
$D({S})$ and right modules, right ideals, left modules, left ideals
of $D(\omega({S}))$, respectively.
\end{theorem}

\begin{proof}
This is immediate from Theorem \ref{theorem:SurjectivityByS_2}.
\end{proof}

\begin{theorem}
\label{thm:LeftNoetherian}
If ${S}$ satisfies condition \eqref{condition:LeftNoether}, then
$D({S})$ is left Noetherian.
In particular, $D({S})$ is left Noetherian,
if ${S}$ satisfies Serre's condition $(S_2)$.
\end{theorem}

\begin{proof}
This is immediate from Corollary \ref{cor:omegaS'isRightNoetherian}
and Theorem \ref{theorem:Left-Right}.
\end{proof}

\begin{theorem}
\label{SelfDuality}
Assume that ${S}$ satisfies condition \eqref{condition:LeftNoether}
and that there exists $\a\in \Z^d$ such that
$\omega({S})=\a+{S}$.
Then
there exist one-to-one correspondences between
left modules, left ideals of $D({S})$ and its
right modules, right ideals, respectively.
\end{theorem}

\begin{proof}
We have
$$
D(\omega({S}))= D(\a+{S})
= t^\a D({S})t^{-\a} \simeq  D({S})
$$
as ${K}$-algebras.
Hence the theorem follows from Theorem \ref{theorem:Left-Right}.
\end{proof}

\begin{corollary}
\label{Gorenstein}
Assume that $K[\N A]$ is Gorenstein.
Then
there exist one-to-one correspondences between
left modules, left ideals of $D(\N A)$ and its
right modules, right ideals, respectively.
\end{corollary}

\begin{proof}
In this case, $\N A$ satisfies the assumption of Theorem \ref{SelfDuality}.\end{proof}

\begin{example}
\label{example:1} 
Let
$$
A = 
\left[
\begin{array}{llllllll} 
0 & 1 & 3 & 3 & 4 & 4 & 6 & 6\\
2 & 1 & 1 & 2 & 0 & 1 & 0 & 1
\end{array} \right]. 
$$
The semigroup $\N A$ is illustrated in Figure \ref{fig2}.
The scored extention $S_c(\N A)$ equals $\N^2$.
We have
$$
\Z(A\cap\sigma_1)
=
\Z\begin{bmatrix}
2 \\ 0
\end{bmatrix}
\quad
\text{and}
\quad
\Z(A\cap\sigma_2)
=
\Z\begin{bmatrix}
0 \\ 2
\end{bmatrix}.
$$
The standard expression of $\N A$ is
$$
\N A =\N^2\setminus
\left(
\begin{bmatrix}
1 \\ 0
\end{bmatrix}
+\Z
\begin{bmatrix}
2 \\ 0
\end{bmatrix}
\right)
\setminus
\left(
\begin{bmatrix}
0 \\ 1
\end{bmatrix}
+\Z
\begin{bmatrix}
0 \\ 2
\end{bmatrix}
\right)
\setminus
\left(
\begin{bmatrix}
1 \\ 0
\end{bmatrix}
+\Z
\begin{bmatrix}
0 \\ 2
\end{bmatrix}
\right)
\setminus
\left(
\begin{bmatrix}
2 \\ 1
\end{bmatrix}
+\Z
\begin{bmatrix}
0 \\ 2
\end{bmatrix}
\right)
\setminus
\begin{bmatrix}
2 \\ 0
\end{bmatrix}.
$$
The $S_2$-closure $S_2(\N A)$ equals $\N A\cup 
{\displaystyle
\left\{
\begin{bmatrix}
2 \\ 0
\end{bmatrix}
\right\}}$, and thus
$\N A$ does not satisfy $(S_2)$ condition.

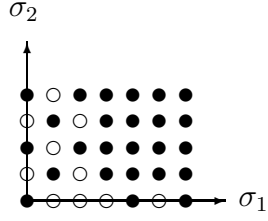
\begin{figure}[!ht]
\begin{picture}(100,90)(0,10)

\put(10,10){\circle*{5}}
\put(20,10){\circle{5}}
\put(30,10){\circle{5}}
\put(40,10){\circle{5}}
\put(50,10){\circle*{5}}
\put(60,10){\circle{5}}
\put(70,10){\circle*{5}}
\put(10,20){\circle{5}}
\put(20,20){\circle*{5}}
\put(30,20){\circle{5}}
\put(40,20){\circle*{5}}
\put(50,20){\circle*{5}}
\put(60,20){\circle*{5}}
\put(70,20){\circle*{5}}
\put(10,30){\circle*{5}}
\put(20,30){\circle{5}}
\put(30,30){\circle*{5}}
\put(40,30){\circle*{5}}
\put(50,30){\circle*{5}}
\put(60,30){\circle*{5}}
\put(70,30){\circle*{5}}
\put(10,40){\circle{5}}
\put(20,40){\circle*{5}}
\put(30,40){\circle{5}}
\put(40,40){\circle*{5}}
\put(50,40){\circle*{5}}
\put(60,40){\circle*{5}}
\put(70,40){\circle*{5}}
\put(10,50){\circle*{5}}
\put(20,50){\circle{5}}
\put(30,50){\circle*{5}}
\put(40,50){\circle*{5}}
\put(50,50){\circle*{5}}
\put(60,50){\circle*{5}}
\put(70,50){\circle*{5}}
\put(10,10){\vector(1,0){75}}
\put(10,10){\vector(0,1){60}}
\put(3,80){$\sigma_2$}
\put(90,8){$\sigma_1$}
\end{picture}
\caption{The semigroup ${\N} A$ in Example \ref{example:1}}
\label{fig2}
\end{figure}
We have
$
\tilde{\cF}=
\{
\sigma_1, \sigma_2, \{\0\}
\}$.\newline
The set $\cB$ consists of only one element $(\b_\tau):$
\begin{picture}(140,30)(0,50)
\put(10,50){$\0=\b_{\{\0\}}$}
\put(100,75){$\b_{\sigma_1}=\begin{bmatrix}
1 \\ 0
\end{bmatrix}$}
\put(100,25){$\b_{\sigma_2}=\begin{bmatrix}
0 \\ 1
\end{bmatrix}$.}
\put(55,55){\line(2,1){40}}
\put(55,55){\line(2,-1){40}}
\end{picture}
\vspace{50pt}

\noindent
We can check that, for this $(\b_\tau)$,
${\displaystyle
\omega(\N A)=\begin{bmatrix}
-2 \\ 0
\end{bmatrix}
+\N A.}
$
Hence, by Theorems \ref{theorem:Left-Right} and \ref{thm:LeftNoetherian},
the category of left $D(\N A)$-modules and that
of right $D(\N A)$-modules are equivalent, and
$D(\N A)$ is Noetherian.
\end{example}

\section{Left Noetherian Property}
\label{LeftNoether}

In this section,
we consider the left Noetherian property;
we give a sufficient condition in Subsection \ref{subsection:SufficientCondition}
and a necessary condition
in Subsection \ref{subsection:NecessaryCondition}.

Let $S$ be a semigroup $\N A$,
and $S_2$ its $S_2$-closure.

Let 
\begin{equation}
\label{StandardExpressionSandS2}
\begin{array}{l}
{\displaystyle S=S_c\setminus \bigcup_{i=1}^m(\b_i+\Z(A\cap\tau_i))}\\
{\displaystyle S_2=S_c\setminus \bigcup_{i=1}^l(\b_i+\Z(A\cap\tau_i))}
\end{array}
\end{equation}
be the standard expressions for $S$ and $S_2$.
Hence $\tau_i$ $(i\leq l)$ are facets, while $\tau_i$ $(l<i\leq m)$ are not.
\begin{lemma}
$$
D(S)\subseteq D(S_2).
$$
\end{lemma}

\begin{proof}
Recall that
\begin{eqnarray*}
D(S_2)_\a &=& 
D(S_c)_\a\cdot\I(\bigcup_{\b_i-\a\in S_2+\Z(A\cap\tau_i),\, i\leq l}
\b_i-\a+{K}\tau_i),\\
D(S)_\a &=& 
D(S_c)_\a\cdot\I(\bigcup_{\b_i-\a\in S+\Z(A\cap\tau_i),\, i\leq m}
\b_i-\a+{K}\tau_i).
\end{eqnarray*}
Note that for $i\leq l$
$$
S_2+\Z(A\cap\tau_i)
=
S+\Z(A\cap\tau_i)
$$
since $\tau_i$ is a facet.
Hence $D(S)\subseteq D(S_2)$.
\end{proof}

\subsection{A Sufficient Condition}
\label{subsection:SufficientCondition}

Since $D(S_2)$ is left Noetherian by Theorem \ref{thm:LeftNoetherian},
the following lemma is proved similarly to Lemma \ref{Robson'sLemma}.

\begin{lemma}
\label{LeftRobson'sLemma}
If $D(S_2)/D(S_2, S)$ is a Noetherian left $D(S)$-module,
then $D(S)$ is left Noetherian.
\end{lemma}

\begin{theorem}
\label{LeftNoetherianFacetAndOrigin}
Assume that $S_2\setminus S$ is a finite set.
If, for all $i>l$, the intersection 
$$
\bigcap_{\b_i-\b_j\in{K}\tau_j,\, j\leq l}\tau_j
$$
equals the origin,
then
$D(S)$ is left Noetherian.
\end{theorem}

\begin{proof}
We show that $D(S_2)/D(S_2,S)$ is finite-dimesional.
Then the theorem follows from Lemma \ref{LeftRobson'sLemma}.

Note that all $\tau_i$ with $i>l$ are the origin $\{ \0\}$,
since $S_2\setminus S$ is finite.

First we show that $D(S_2)_\a=D(S_2,S)_\a$ for all but finite $\a\in\Z^d$.
Recall that
\begin{eqnarray*}
D(S_2)_\a &=& 
D(S_c)_\a\cdot\I(\bigcup_{\b_i-\a\in S_2+\Z(A\cap\tau_i),\, i\leq l}
\b_i-\a+{K}\tau_i),\\
D(S_2, S)_\a &=& 
D(S_c)_\a\cdot\I(\bigcup_{\b_i-\a\in S_2+\Z(A\cap\tau_i),\, i\leq m}
\b_i-\a+{K}\tau_i).
\end{eqnarray*}
Hence $D(S_2)_\a \neq D(S_2,S)_\a$ if and only if
there exists $i>l$ such that
\begin{equation}\label{Weights}
\begin{array}{l}
\bullet\quad \b_i-\a\in S_2,\\
\bullet\quad \b_i-\b_j\in{K}\tau_j\, (j\leq l)
\Rightarrow \b_j-\a\notin S_2+\Z(A\cap\tau_j).
\end{array}
\end{equation}
It suffices to show that for a fixed $i>l$ there exists only finitely many
$\a\in \Z^d$ with \eqref{Weights}.
Take $M$ as in \eqref{Def-M}, i.e.,
$$
M> \max F_\sigma(S)^c\cup\{ F_\sigma(\b_i)\, :\, i\}
\quad\text{for all facets $\sigma$}.
$$
(Note that
$\min F_\sigma(S)\cup\{ F_\sigma(\b_i)\, :\, i\}=0$
in our case.)

Suppose that $\a\in \Z^d$ satisfies \eqref{Weights}.
Then
\begin{equation}\label{Weights2}
\begin{array}{l}
\bullet\quad F_\sigma(\b_i-\a)\geq 0\quad\text{for all facets $\sigma$,}\\
\bullet\quad F_{\tau_j}(\b_i-\a)=F_{\tau_j}(\b_j-\a)
\leq M\quad\text{for all $j\leq l$ with
$\b_i-\b_j\in{K}\tau_j$.}
\end{array}
\end{equation}
There exists only finitely many such $\a\in\Z^d$, since the intersection
$
{\displaystyle
\bigcap_{\b_i-\b_j\in{K}\tau_j,\, j\leq l}\tau_j}
$
equals the origin.

Now it is left to show that
each $D(S_2)_\a/D(S_2, S)_\a$ is finite-dimensional. 
Let
$I:=\I(\{ \b_i-\a\,:\,
\b_i-\a\in S_2,\, l< i\leq m\})$.
Then $D(S_2)_\a I\subseteq D(S_2,S)_\a$.
There exist surjective ${K}[s]$-module homomorphisms
$$
D(S_2)_\a/D(S_2,S)_\a
\leftarrow D(S_2)_\a/D(S_2)_\a I
\leftarrow {K}[s]/I.
$$
The latter is an isomorphism, since $D(S_2)_\a$ is a singly generated
${K}[s]$-module by \cite[Proposition 7.7]{Prim}.
Hence $D(S_2)_\a/D(S_2, S)_\a$ is finite-dimensional,
and we have completed the proof.
\end{proof}
\subsection{A Necessary Condition}
\label{subsection:NecessaryCondition}

Let $S$ be a semigroup $\N A$, and let
$$
{S}= S_c\setminus \bigcup_{i=1}^m(\b_i+\Z(A\cap\tau_i))
$$
be the standard expression, where $S_c$ is the scored extention of ${S}$.
In this subsection,
we assume that $\codim\, \tau_m>1$, and that
\begin{equation}
\label{Condition:nonLeftNoether}
\left(\bigcap_{\tau_i\succ \tau_m;\, \b_i-\b_m\in{K}\tau_i}\tau_i\right)
\neq\tau_m,
\end{equation}
and we show that $D(S)$ is not left Noetherian.
We construct a strictly increasing sequence of left ideals of $D(S)$.

Let $\rho$ be a ray of $\R_{\geq 0}A$ contained in $(\bigcap_{\tau_i\succ \tau_m;\, \b_i-\b_m\in{K}\tau_i}\tau_i)
\setminus(\tau_m\setminus \{\0\})$.
Fix a vector $\d_\rho\in \N(A\cap\rho)$.
Similarly to Lemma \ref{lemma:Mk}, we have the following Lemma.

\begin{lemma}
\label{lemma:m-rho}
\begin{equation}
D({S})_{-{k}\d_{\rho}}=
\I(\bigcup_{\b_i+{k}\d_\rho\in{S}+\Z(A\cap\tau_i)}(\b_i+{K}\tau_i))
\cdot P_{-{k}\d_\rho},
\end{equation}
where $D(S_c)_{-{k}\d_{\rho}}={K}[s] P_{-{k}\d_\rho}$.
\end{lemma}

\begin{lemma}
\label{lemma:1}
For ${k}\gg 0$,
\begin{equation}
\b_m+{k}\d_\rho\in{S}+\Z(A\cap\tau_m).
\end{equation}
\end{lemma}

\begin{proof}
Since $\d_\rho\in {S}$,
$\b_m+{k}\d_\rho\in S_c$.
Suppose that
$\b_m+{k}\d_\rho\not\in {S} +\Z(A\cap\tau_m)$.
Then there exists $i$ with $\tau_i\succeq\tau_m$
such that $\b_m+{k}\d_\rho\in\b_i+\Z(A\cap\tau_i)$
for ${k}\gg 0$.
Hence $\d_\rho\in\Z(A\cap\tau_i)$.
Thus we have $\b_m\in\b_i+\Z(A\cap\tau_i)$,
and then $\b_m+\Z(A\cap\tau_m)\subseteq
\b_i+\Z(A\cap\tau_i)$.
By the irredundancy of the standard expression, we have $i=m$.
But, since $\rho\not\preceq \tau_m$, we have
$\b_m+{k}\d_\rho\not\in \b_m+\Z(A\cap\tau_m)$.
\end{proof}

\begin{lemma}
\label{lemma:i}
Suppose that
$\tau_i\succ\tau_m$ and $\b_i-\b_m\in {K}\tau_i$.
Then
\begin{equation}
\b_i+{k}\d_\rho\not\in{S}+\Z(A\cap\tau_i).
\end{equation}
\end{lemma}

\begin{proof}
By the definition of $\rho$,
$\rho\preceq \tau_i$ for such $\tau_i$.
Hence $\d_\rho\in \Z(A\cap\tau_i)$.
Then we see $\b_i+{k}\d_\rho\not\in {S}+\Z(A\cap\tau_i)$.
\end{proof}

For each $i$ with $\b_i+{k}\d_\rho\in{S}+\Z(A\cap\tau_i)$
and $\b_i+{K}\tau_i\neq \b_m+{K}\tau_m$,
we take a facet $\sigma_i\succeq \tau_i$ as follows:
\begin{enumerate}
\item
If $\tau_i\not\succeq \tau_m$, then take a facet $\sigma_i\succeq\tau_i$
such that $\sigma_i\not\succeq\tau_m$.
\item
If $\tau_i\succeq\tau_m$ and $\b_i-\b_m\not\in{K}\tau_i$, then
take a facet $\sigma_i\succeq\tau_i$ such that $F_{\sigma_i}(\b_i)\neq
F_{\sigma_i}(\b_m)$.
\item
We do not need to consider the case where
$\tau_i\succ\tau_m$ and $\b_i-\b_m\in{K}\tau_i$ by Lemma \ref{lemma:i}.
\end{enumerate}
Finally take a facet $\sigma_m$ containing $\rho$ and $\tau_m$.
Then define $E({k})$ by
\begin{equation}
\label{equation:E({m})}
E({k}):=
(F_{\sigma_m}-F_{\sigma_m}(\b_m))
\prod_{\overset{\scriptstyle\b_i+{k}\d_\rho\in{S}+\Z(A\cap\tau_i),}
{\b_i+{K}\tau_i\neq \b_m+{K}\tau_m}}
(F_{\sigma_i}-F_{\sigma_i}(\b_i))\cdot
P_{-{k}\d_\rho}.
\end{equation}

Then, by Lemma \ref{lemma:m-rho}, $E({k})\in D({S})_{-{k}\d_\rho}$.

\begin{lemma}
\label{lemma:E(m)notin}
$E({k})\not\in (F_{\sigma_m}-F_{\sigma_m}(\b_m)) D({S})_{-{k}\d_\rho}$
for $k\gg 0$.
\end{lemma}

\begin{proof}
By the definitions of $\sigma_i$,
$\prod_{\overset{\scriptstyle\b_i+{k}\d_\rho\in{S}+\Z(A\cap\tau_i),}
{\b_i+{K}\tau_i\neq \b_m+{K}\tau_m}}
(F_{\sigma_i}-F_{\sigma_i}(\b_i))
\not\in
\I(\b_m+{K}\tau_m)$.
Hence the assertion follows from Lemmas \ref{lemma:m-rho} and \ref{lemma:1}.
\end{proof}

\begin{theorem}
\label{NecessaryCondition}
$D({S})$ is not left Noetherian.
\end{theorem}

\begin{proof}
We construct a strictly increasing sequence of left ideals.

First we claim that for $k,l\gg 0$
\begin{equation}
\label{claim}
D({S})_{-l\d_\rho}\cdot E({k})
\subseteq
(F_{\sigma_m}-F_{\sigma_m}(\b_m)) D({S})_{-(l+{k})\d_\rho}.
\end{equation}
Let $f(s)P_{-l\d_\rho}\in D({S})_{-l\d_\rho}$.
Note that, by Lemma \ref{lemma:m-rho}, for $k,l\gg 0$ we have
$f(s)P_{-k\d_\rho}\in D({S})_{-k\d_\rho}$
if and only if $f(s)P_{-l\d_\rho}\in D({S})_{-l\d_\rho}$,
and recall from \cite[Lemma 5.10]{Saito-Traves2} that
$P_{-(l+k)\d_\rho}=P_{-l\d_\rho}P_{-k\d_\rho}$ for $k,l\gg 0$.
Thus for $k,l\gg 0$
\begin{eqnarray*}
f(s)P_{-l\d_\rho}E({k})
&=&
f(s)P_{-l\d_\rho}(F_{\sigma_m}-F_{\sigma_m}(\b_m))
\prod_{\overset{\scriptstyle\b_i+{k}\d_\rho\in{S}+\Z(A\cap\tau_i),}
{\b_i+{K}\tau_i\neq \b_m+{K}\tau_m}}
(F_{\sigma_i}-F_{\sigma_i}(\b_i))
P_{-{k}\d_\rho}\\
&\in&
(F_{\sigma_m}-F_{\sigma_m}(\b_m))
{K}[s]
f(s)P_{-l\d_\rho}P_{-{k}\d_\rho}\\
&=&
(F_{\sigma_m}-F_{\sigma_m}(\b_m))
{K}[s]
f(s)P_{-(l+{k})\d_\rho}\\
&\subseteq&
(F_{\sigma_m}-F_{\sigma_m}(\b_m)) D({S})_{-(l+{k})\d_\rho}.
\end{eqnarray*}
Thus we have proved claim \eqref{claim}.

Take $k_0\in\N$ large enough.
By Lemma \ref{lemma:E(m)notin} and claim \eqref{claim},
for $l>h$,
$$
E(lk_0)\notin (\sum_{k=1}^hD({S})\cdot E({kk_0}))_{-lk_0\d_\rho}.
$$
Therefore
$\{ \sum_{{k}=1}^hD({S})\cdot E({kk_0})\,:\, h=1,2,\ldots\}$
is a strictly increasing sequence of left ideals of $D({S})$.
\end{proof}

Finally we make a conjecture of the condition
for $D(S)$ to be left Noetherian.

\begin{conjecture}
Let \eqref{StandardExpressionSandS2} be 
the standard expressions of $S$ and $S_2$.
Then the following are equivalent.
\begin{enumerate}
\item[{\rm (1)}]
$D(S)$ is left Noetherian.
\item[{\rm (2)}]
$D(S_2)/D(S_2,S)$ is a Noetherian left $D(S)$-module.
\item[{\rm (3)}]
$D(S_2)/D(S)$ is a Noetherian left $D(S)$-module.
\item[{\rm (4)}]
For all $i>l$, 
$$
\left(
\bigcap_{\tau_j\succ\tau_i,\, \b_i-\b_j\in {K}\tau_j,\,
j\leq l}\tau_j
\right)
=\tau_i.
$$
\end{enumerate}
\end{conjecture}

\begin{remark}
Lemma \ref{LeftRobson'sLemma} says that $(2)$ implies $(1)$.
Clearly $(2)$ and $(3)$ are equivalent under $(1)$.
The implication $(1)\Rightarrow (4)$ is Theorem \ref{NecessaryCondition}.
Note also that $(4)$ is satisfied when the set $\cB$ is not empty
(cf. Theorem \ref{thm:LeftNoetherian}).

By Theorems \ref{LeftNoetherianFacetAndOrigin} and 
\ref{NecessaryCondition}, the conjecture is true for $d\leq 2$.
\end{remark}

 
\ifx\undefined\bysame 
\newcommand{\bysame}{\leavevmode\hbox to3em{\hrulefill}\,} 
\fi

\end{document}